\documentclass{article}

\usepackage{amssymb,exscale,amsmath}

\author{Steffen Brasch\footnote{Department of Mathematics and Computer
    Science, University of Greifswald, Jahnstr. 15a,
    17487 Greifswald, Germany,
    {\it e-mail: sbrasch@uni-greifswald.de}},
  Katharina Habermann\footnote{SUB G\"ottingen,
    Platz der G\"ottinger Sieben 1, 37073 G\"ottingen, Germany,
    {\it e-mail: habermann@sub.uni-goettingen.de}},
  Lutz Habermann\footnote{Institute of Differential Geometry,
    University of Hannover, Welfengarten 1, 30167 Hannover, Germany,
    {\it e-mail: habermann@math.uni-hannover.de}}}

\title{Symplectic Dirac Operators on Hermitian Symmetric
  Spaces\footnote{This work was partially supported by the
    Deutsche Forschungsgemeinschaft (DFG).}}

\date{\empty}

\numberwithin{equation}{section}

\newcommand{\vp}{\varphi}
\newcommand{\rd}{\mathrm{d}}
\newcommand{\C}{\mathbb{C}}
\newcommand{\N}{\mathbb{N}}
\newcommand{\R}{\mathbb{R}}
\newcommand{\mfg}{\mathfrak{g}}
\newcommand{\mfh}{\mathfrak{h}}
\newcommand{\mfm}{\mathfrak{m}}
\newcommand{\mfsu}{\mathfrak{su}}
\newcommand{\mfu}{\mathfrak{u}}
\newcommand{\mfV}{\mathfrak{V}}
\newcommand{\mfW}{\mathfrak{W}}
\newcommand{\bfU}{\mathbf{U}}
\newcommand{\sfa}{\mathsf{a}}
\newcommand{\sfv}{\mathsf{v}}
\newcommand{\sfA}{\mathsf{A}}
\newcommand{\re}{\mathrm{e}}
\newcommand{\ri}{\mathrm{i}}
\newcommand{\GL}{\mathrm{GL}}
\newcommand{\Hom}{\mathrm{Hom}}
\newcommand{\Mp}{\mathrm{Mp}}
\newcommand{\rT}{\mathrm{T}}
\newcommand{\Sp}{\mathrm{Sp}}
\newcommand{\SU}{\mathrm{SU}}
\newcommand{\rU}{\mathrm{U}}
\newcommand{\Sec}{{\mathrm{\Gamma}}}
\newcommand{\mcD}{\mathcal{D}}
\newcommand{\mcP}{\mathcal{P}}
\newcommand{\Ad}{\operatorname{Ad}}
\newcommand{\ad}{\operatorname{ad}}

\newtheorem{lm}{Lemma}[section]
\newtheorem{sz}[lm]{Proposition}
\newtheorem{tm}[lm]{Theorem}
\newtheorem{fo}[lm]{Corollary}

\newtheorem{bm0}[lm]{Remark}
\newtheorem{bs0}[lm]{Example}

\newenvironment{beweis}{{\em Proof.}}{\hspace*{\fill} $\square$}

\newenvironment{bs}{\begin{bs0}\rm}{\hspace*{\fill} $\square$ \end{bs0}}

\begin{document}

\maketitle

\begin{abstract}
  We describe the shape of the symplectic Dirac operators on Hermitian
  symmetric spaces. For this, we consider these operators as
  families of operators that can be handled more easily than the original ones.
\end{abstract}

\section{Introduction}

Symplectic spinor fields were introduced by B.~Kostant in the early 1970th
\cite{kostant}.
This was done in the context of geometric quantization, where he was
interested in constructing the so-called half-form bundle as well as the
half-form pairing in order to establish the appropriate Hilbert space.
Kostants construction was quite an essential step in geometric quantization,
but afterwards no further work has been done to seek
out interior aspects of symplectic spin geometry for two decades.

In 1995, the second author introduced symplectic Dirac operators and
initialized a systematical investigation \cite{khab1}.
Some progress has been made in understanding basic properties of these
operators \cite{khab2}.
Moreover, symplectic Dirac operators have been studied in the setting of
parabolic geometries (cf. e.g. \cite{lenka,krysl1}) and received
attention in mathematical physics (cf. e.g. \cite{sommen}).
A systematic and thorough introduction to the subject can be found in
\cite{hab2-buch}.

There are several topics of their own interest that are relevant to get a
broader understanding of symplectic Dirac operators.
For example, a deeper knowledge of the moduli space of symplectic structures
may provide a way of analyzing how the symplectic Dirac operators
are related to the choice of the symplectic structure of the underlying
symplectic manifold. For considerations on moduli spaces of symplectic
structures, see \cite{hab1-jan,hab2-jan,wilson}.
Furthermore, the symplectic Dirac operators are defined by means of a
symplectic connection.
It is well known that in symplectic geometry there is no analog of the
Levi-Civita connection.
With the aim to get a deeper insight into the structure of the space of
symplectic connections, in \cite{hab2-paul} an approach of a purely symplectic
Yang-Mills theory is given.

In order to develop the symplectic Dirac operator picture, it is
significant to have more explicit examples.
In the present paper, we study symplectic Dirac operators in the case where
the underlying symplectic manifold is a Hermitian symmetric space.
This gives a class of examples which were not figured out so far.
In the particular case of odd-dimensional complex projective spaces, the
spectrum of an associated second order operator has been already computed
\cite{khab,wyss}.
However, no calculations are performed for the first order symplectic Dirac
operators there.
Here, we are going to study the shape of these operators itself.
Our strategy is based on the idea of considering symplectic Dirac operators
as families of operators acting on more suitable spaces than the section space
of an infinite rank vector bundle.

The paper is organized as follows.
The necessary concepts of symplectic Dirac operators in the setting of
symmetric symplectic spaces are introduced in Section~\ref{sec:prelim}.
For the general situation, we refer to \cite{hab2-buch}.
In Section~\ref{sec:decom} we describe a splitting of the
symplectic spinor bundle into subbundles of finite rank and apply
the Frobenius reciprocity to decompose the
corresponding section spaces in the case of a Hermitian symmetric space.
Moreover, we explain how the symplectic Dirac operators built with respect
to the canonical Hermitian connection behave in relation to
these decompositions. The last section is devoted to a detailed discussion
of the case $\C P^1$.

\section{Preliminaries
\label{sec:prelim}}

Let $G$ be a simply connected real Lie group, let $H$ be a closed and
connected subgroup of $G$ and set $M$ to be the homogeneous space $G/H$.
Let $\mfg$ and $\mfh$ denote the Lie algebras of $G$ and $H$, respectively.
We suppose that there exists a subspace $\mfm\subset\mfg$ such that
\begin{equation}
  \label{eq:zerleg-lag}
  \mfg=\mfh\oplus\mfm
\end{equation}
as well as
\begin{displaymath}
  [\mfh,\mfm]\subset\mfm
  \quad\mbox{and}\quad
  [\mfm,\mfm]\subset\mfh\;.
\end{displaymath}
Moreover, let $\omega_0$ be an $\ad(\mfh)$-invariant non-degenerate
skew-symmetric bilinear form on $\mfm$. Identifying the tangent
space $T_o M$ to $M$ at the point $o=eH$ with $\mfm$, $\omega_0$ induces
a $G$-invariant almost symplectic structure $\omega$ on $M$.
Here and in the sequel, $e$ denotes the unit element of $G$.
Since the $\ad(\mfh)$-invariance of $\omega_0$ is equivalent
to the trivial extension of $\omega_0$ to $\mfg$ being
a Chevalley $2$-cocycle, $\omega$ is closed. Thus $\omega$ is even a
symplectic structure on $M$. In this way, we have described
$(M,\omega)$ as a simply connected symplectic symmetric space
(cf. \cite{bieliavsky}).

The canonical projection $\pi:G\to M$ and the right action of $H$ on $G$
give $G$ the structure of an $H$-principal fiber bundle.
Let $\kappa:H\to\GL(\mfm)$ be the isotropy representation of the
homogeneous space $M$, i.e. the restriction of the adjoint representation
of $G$ to $H$ acting on $\mfm$. By assumption, $\kappa$ maps into the
symplectic group $\Sp(\mfm)$ of the symplectic vector
space $(\mfm,\omega_0)$. We fix a symplectic basis $(X_1,\ldots,X_{2n})$
of $(\mfm,\omega_0)$, i.e. a basis of $\mfm$ such that
\begin{displaymath}
  \omega_0(X_j,X_k)=\omega_0(X_{n+j},X_{n+k})=0
  \quad\mbox{and}\quad
  \omega_0(X_j,X_{n+k})=\delta_{jk}
\end{displaymath}
for $j,k=1,\ldots,n$, and identify
the symplectic frame bundle $R$ of $(M,\omega)$ with the
$\Sp(\mfm)$-principal fiber bundle $G\times_\kappa\Sp(\mfm)$ associated
to $G$ with respect to $\kappa:H\to\Sp(\mfm)$ via
\begin{displaymath}
  [g,A]\in G\times_\kappa\Sp(\mfm)
  \mapsto(\rd\pi(\rd L_g (AX_1)),\ldots,\rd\pi(\rd L_g (AX_{2n})))\in R\;.
\end{displaymath}
Let $\Mp(\mfm)$ be the connected double covering group of $\Sp(\mfm)$
and let $\rho:\Mp(\mfm)\to\Sp(\mfm)$ denote the covering homomorphism.
If $\tilde\kappa:H\to\Mp(\mfm)$ is a lift of $\kappa$, i.e.
a homomorphism such that $\rho\circ\tilde\kappa=\kappa$, then the
$\Mp(\mfm)$-principal fiber bundle $P=G\times_{\tilde\kappa}\Mp(\mfm)$
together with the map $F_P:P\to R$ defined by
$F_P([g,q])=[g,\rho(q)]$ is a metaplectic structure of the symplectic
manifold $(M,\omega)$.

\begin{tm}
  The mapping $\tilde\kappa\mapsto (P,F_P)$ described above induces
  a $1:1$ correspondence between the lifts $\tilde\kappa$ of $\kappa$
  and the isomorphism classes of metaplectic structures $(P,F_P)$
  of the simply connected symplectic symmetric space $(M,\omega)$.
\end{tm}
\begin{beweis}
  One proceeds as in the Riemannian case (cf. \cite{baer,cah-gu,steffen}).
\end{beweis}

\medskip

According to the above, we identify the tangent bundle $TM$ of $M$
with the vector bundle $G\times_\kappa\mfm$ associated to $G$ with respect
to $\kappa$ via
\begin{displaymath}
  [g,X]\in G\times_\kappa\mfm\mapsto\rd\pi(\rd L_g X)\in TM\;.
\end{displaymath}
The space of smooth vector fields on $M$, i.e. the space
$\Sec(TM)$ of smooth sections of $TM$ is then the space of
smooth maps $\xi:G\to\mfm$ such that
\begin{displaymath}
  \xi(gh)=\kappa\left(h^{-1}\right)\xi(g)
\end{displaymath}
for $g\in G$ and $h\in H$ and the symplectic structure $\omega$ is given by
\begin{displaymath}
  \omega([g,X],[g,Y])=\omega_0(X,Y)\;.
\end{displaymath}
Consequently, the function $\omega(\xi_1,\xi_2)\in C^\infty(M)$ for
$\xi_1,\xi_2\in\Sec(TM)$ satisfies
\begin{equation}
  \label{eq:rel-om-Om}
  \omega(\xi_1,\xi_2)(\pi(g))=\omega_0(\xi_1(g),\xi_2(g))\;.
\end{equation}

Let $\nabla$ be the $G$-invariant connection on $M$ defined by
\begin{displaymath}
  \nabla_{[g,X]}\xi=[g,\rd\xi(\rd L_g X)]
\end{displaymath}
for $[g,X]\in TM$ and $\xi\in\Sec(TM)$. Then
\begin{equation}
  \label{eq:descr-sympl-conn}
  (\nabla_{\xi_1}\xi_2)(g)=\rd\xi_2(\rd L_g\xi_1(g))
\end{equation}
for $\xi_1,\xi_2\in\Sec(TM)$ and $g\in G$.

\begin{sz}
  The connection $\nabla$ is symplectic, i.e. $\nabla\omega=0$.
\end{sz}
\begin{beweis}
  Let $\xi,\xi_1,\xi_2\in\Sec(TM)$. By means of
  \begin{displaymath}
    \xi(u)(\pi(g))=\rd(u\circ\pi)(\rd L_g\xi(g))
  \end{displaymath}
  for $u\in C^\infty(M)$ and Equations~(\ref{eq:rel-om-Om}) and
  (\ref{eq:descr-sympl-conn}), we conclude
  \begin{align*}
    \xi(\omega(\xi_1,\xi_2))(\pi(g))
    &=\rd(\omega_0(\xi_1,\xi_2))(\rd L_g\xi(g))\\
    &=\omega_0(\rd\xi_1(\rd L_g\xi(g)),\xi_2(g))
    +\omega_0(\xi_1(g),\rd\xi_2(\rd L_g\xi(g)))\\
    &=\omega_0((\nabla_\xi\xi_1)(g),\xi_2(g))
    +\omega_0(\xi_1(g),(\nabla_\xi\xi_2)(g))\\
    &=\omega(\nabla_\xi\xi_1,\xi_2)(\pi(g))
    +\omega(\xi_1,\nabla_\xi\xi_2)(\pi(g))\;.
  \end{align*}
\end{beweis}

\medskip

Moreover, we have (cf. \cite{ko-no2,loos})

\begin{sz}
  The connection $\nabla$ is torsion-free.
  \hspace*{\fill} $\square$
\end{sz}

Thus $\nabla$ is a torsion-free symplectic connection on the almost
symplectic manifold $(M,\omega)$. By a classical result of Tondeur
(cf. \cite{tondeur}), this also implies that $(M,\omega)$ is symplectic.

Using the fixed symplectic basis $(X_1,\ldots,X_{2n})$ of $(\mfm,\omega_0)$,
we identify $\Sp(\mfm)$ and its double cover $\Mp(\mfm)$
with the symplectic group $\Sp(n,\R)$ and the metaplectic group
$\Mp(n,\R)$, respectively. Let
$\boldsymbol{m}:\Mp(\mfm)\to\rU(L^2(\R^n))$ be the metaplectic representation.
Let $\tilde\kappa$ be a lift of $\kappa$ and let $(P,F_P)$
be the metaplectic structure of $(M,\omega)$ constructed from $\tilde\kappa$
as described above. The symplectic spinor bundle is then the
Hilbert space bundle $Q=G\times_\lambda L^2(\R^n)$ associated to
$G$ with respect to
$\lambda=\boldsymbol{m}\circ\tilde\kappa:H\to\rU(L^2(\R^n))$.
Accordingly, a symplectic spinor field, i.e. a smooth section of $Q$
is a smooth map $\vp:G\to L^2(\R^n)$ such that
\begin{displaymath}
  \vp(gh)=\lambda\left(h^{-1}\right)\vp(g)
\end{displaymath}
for $g\in G$ and $h\in H$. Since the symplectic connection $\nabla$
on $M$ is induced
by the $G$-invariant connection on the $H$-principal fiber bundle $G$
that corresponds to the decomposition (\ref{eq:zerleg-lag}),
the same holds true for the spinor derivative on $Q$, also denoted by
$\nabla$. Therefore,
\begin{equation}
  \label{eq:express-spinor-der}
  (\nabla_\xi\vp)(g)=\rd\vp(\rd L_g\xi(g))
\end{equation}
for $\vp\in\Sec(Q)$ and $\xi\in\Sec(TM)$.

According to our conventions, the symplectic Clifford multiplication
is the homomorphism $\mu_0:\mfm\otimes L^2(\R^n)\to L^2(\R^n)$
given by
\begin{equation}
  \label{eq:def-cliff-mult}
  \mu_0(X_j\otimes f)(x)=\ri x_j f(x)
  \quad\mbox{and}\quad
  \mu_0(X_{n+j}\otimes f)
  =\frac{\partial f}{\partial x^j}
\end{equation}
for $f\in L^2(\R^n)$, $x=(x_1,\ldots,x_n)\in\R^n$ and $j=1,\ldots,n$,
where the multiplication by $X\in\mfm$ has to be seen as an unbounded
operator on $L^2(\R^n)$. As usual, we shall write $X\cdot f$ for
$\mu_0(X\otimes f)$. The multiplication $\mu_0$ induces a
multiplication $\mu:TM\otimes Q\to Q$ by
\begin{displaymath}
  \mu([g,X]\otimes[g,f])=[g,X\cdot f]\;.
\end{displaymath}
The first symplectic Dirac operator is now defined as
\begin{displaymath}
  \mcD=\mu\circ(\bar\omega\otimes\operatorname{id}_Q)\circ\nabla:
  \Sec(Q)\to\Sec(Q)\;,
\end{displaymath}
where $\bar\omega$ means the isomorphism $\bar\omega:T^\ast M\to TM$
generated by the symplectic structure $\omega$ on $M$.

\begin{sz}
  \label{sz:D-on-ssymsp}
  The symplectic Dirac operator $\mcD$ can be written as
  \begin{equation}
    \label{eq:express-D}
    \mcD(\vp)=\sum_{j=1}^n X_j\cdot X_{n+j}(\vp)
    -\sum_{j=1}^n X_{n+j}\cdot X_j(\vp)
  \end{equation}
  for $\vp\in\Sec(Q)$, where $X(\vp)$ for $X\in\mfm$ denotes the
  derivative of $\vp$ in the direction of the left-invariant vector
  field determined by $X$, i.e.
  \begin{displaymath}
    X(\vp)(g)=\rd\vp(\rd L_g X)\;.
  \end{displaymath}
\end{sz}
\begin{beweis}
  Let $\xi_1,\ldots,\xi_{2n}\in\Sec(TM)$ such that
  $\xi_j(e)=X_j$ for $j=1,\ldots,2n$. By the local expression for $D$
  (cf. \cite{hab2-buch}) and Equation~(\ref{eq:express-spinor-der}),
  we then have
  \begin{align*}
    \mcD(\vp)(e)
    &=\sum_{j=1}^n\xi_j(e)\cdot\left(\nabla_{\xi_{n+j}}\vp\right)(e)
    -\sum_{j=1}^n\xi_{n+j}(e)\cdot\left(\nabla_{\xi_j}\vp\right)(e)\\
    &=\sum_{j=1}^n X_j\cdot\rd\vp(X_{n+j})
    -\sum_{j=1}^n X_{n+j}\cdot\rd\vp(X_j)\\
    &=\sum_{j=1}^n X_j\cdot X_{n+j}(\vp)(e)
    -\sum_{j=1}^n X_{n+j}\cdot X_j(\vp)(e)\;.
  \end{align*}
  Since both sides of Equation~(\ref{eq:express-D}) describe $G$-invariant
  operators, this proves the proposition.
\end{beweis}

\medskip

To construct the second Dirac operator $\tilde\mcD$, we need an
$\omega$-compatible almost complex structure $J$ on $M$. We assume that
$J$ is $G$-invariant. That is, we suppose that the endomorphism
$J:TM\to TM$ is given by
\begin{displaymath}
  J[g,X]=[g,J_0 X]
\end{displaymath}
for $[g,X]\in TM$, where $J_0$ is an $\ad(\mfh)$-invariant
complex structure on the vector space $\mfm$ such that
\begin{displaymath}
  \mathbf{g}_0(X,Y)=\omega_0(X,J_0 Y)
\end{displaymath}
for $X,Y\in\mfm$ defines an inner product $\mathbf{g}_0$ on $\mfm$.
Then
\begin{equation}
  \label{eq:om-mfm-J-inv}
  \omega_0(J_0 X,J_0 Y)=\omega_0(X,Y)
  \quad\mbox{and}\quad
  \mathbf{g}_0(J_0 X,J_0 Y)=\mathbf{g}_0(X,Y)
\end{equation}
for any $X,Y\in\mfm$. Let $\mathbf{g}$ denote the $G$-invariant
Riemannian metric on $M$ induced by $\mathbf{g}_0$. This means that
\begin{displaymath}
  \mathbf{g}([g,X],[g,Y])=\mathbf{g}_0(X,Y)\;.
\end{displaymath}
By Equation~(\ref{eq:om-mfm-J-inv}),
\begin{displaymath}
  \omega(J\xi_1,J\xi_2)=\omega(\xi_1,\xi_2)
  \quad\mbox{and}\quad
  \mathbf{g}(J\xi_1,J\xi_2)=\mathbf{g}(\xi_1,\xi_2)
\end{displaymath}
for all $\xi_1,\xi_2\in\Sec(TM)$. Furthermore, it follows that $(M,g,J)$
is a Hermitian symmetric space and that $\nabla$ is its
Levi-Civita connection.

The second symplectic Dirac operator is defined as
\begin{displaymath}
  \tilde\mcD
  =\mu\circ(\bar{\mathbf{g}}\otimes\operatorname{id}_Q)\circ\nabla:
  \Sec(Q)\to\Sec(Q)\;,
\end{displaymath}
where $\bar{\mathbf{g}}:T^\ast M\to TM$ is the identification by means of
the Riemannian metric $\mathbf{g}$ on $M$. Analogously to the proof of
Proposition~\ref{sz:D-on-ssymsp}, one can show

\begin{sz}
  \label{sz:tildeD-on-ssymsp}
  Suppose that the basis $(X_1,\ldots,X_{2n})$ of $\mfm$ is unitary, i.e.,
  in addition to symplecticity, it satisfies $J_0 X_j=X_{n+j}$ for
  $j=1,\ldots,n$. Then the symplectic Dirac operator $\tilde\mcD$ takes
  the form
  \begin{displaymath}
    \tilde\mcD(\vp)=\sum_{j=1}^{2n}X_j\cdot X_j(\vp)
  \end{displaymath}
  for $\vp\in\Sec(Q)$.
  \hspace*{\fill} $\square$
\end{sz}

\section{Decompositions and invariant subspaces
\label{sec:decom}}

First we decompose the symplectic spinor bundle $Q$ into subbundles
of finite rank. For this, let $\rU(\mfm)$ denote the unitary group
of the Hermitian vector space $(\mfm,\mathbf{g}_0,J_0)$ and
set $\hat\rU(\mfm)=\rho^{-1}(\rU(\mfm))$,
which is a connected double cover of $\rU(\mfm)$.
The irreducible components of the restriction
$\boldsymbol{u}:\hat\rU(\mfm)\to\rU(L^2(\R^n))$ of the metaplectic
representation $\boldsymbol{m}$ to $\hat\rU(\mfm)$ can be described as
follows. Let $\N_0$ denote the set of non-negative
integers and let $\boldsymbol{h}_\alpha\in L^2(\R^n)$
for a multi-index $\alpha=(\alpha_1,\ldots,\alpha_n)\in\N_0^n$ 
be the Hermite function on $\R^n$ defined by
\begin{equation}
  \label{eq:def-herm-func}
  \boldsymbol{h}_\alpha(x)
  =\boldsymbol{h}_{\alpha_1}(x_1)\cdots\boldsymbol{h}_{\alpha_n}(x_n)\;.
\end{equation}
Here, $\boldsymbol{h}_l$ for $l\in\N_0$ are the classical Hermite functions
on $\R$, which are given by
\begin{displaymath}
  \boldsymbol{h}_l(t)
  =\re^{t^2/2}\frac{\rd^l}{\rd t^l}\left(\re^{-t^2}\right)\;.
\end{displaymath}
As is well known, the Hermite functions $\boldsymbol{h}_\alpha$ form a
complete orthogonal system in $L^2(\R^n)$. Furthermore, for any $l\in\N_0$,
the span $\mfW_l$ of the functions $\boldsymbol{h}_\alpha$ with
$\alpha_1+\cdots+\alpha_n=l$ is an irreducible
$\boldsymbol{u}$-invariant subspace of $L^2(\R^n)$ (cf. \cite{bo-wa}).

Since the complex structure $J_0$ on $\mfm$ is assumed to be
$\ad(\mfh)$-invariant, the homomorphism $\kappa$ maps into $\rU(\mfm)$
and, therefore, the lift $\tilde\kappa$ maps into $\hat\rU(\mfm)$.
Hence the unitary representation $\lambda:H\to\rU(L^2(\R^n))$ can
be written as $\lambda=\boldsymbol{u}\circ\tilde\kappa$. Let
$\boldsymbol{u}_l:\hat\rU(\mfm)\to\rU(\mfW_l)$ be the restriction of
$\boldsymbol{u}$ to the subspace $\mfW_l$ and set
$\lambda_l=\boldsymbol{u}_l\circ\tilde\kappa$.

The above yields

\begin{sz}
  \label{sz:splitting-Q}
  The symplectic spinor bundle $Q$ splits into the orthogonal sum
  of the finite rank subbundles $Q_l=G\times_{\lambda_l}\mfW_l$, $l\in\N_0$.
  \hspace*{\fill} $\square$
\end{sz}

In the following considerations, we want to make use of the Frobenius
reciprocity. To be able to do this, we assume
from now on that $G$ is compact. Let
$\tau_\sfa:G\to\GL(\mfV_\sfa)$, $\sfa\in\sfA$, form a complete system of
representatives of isomorphism classes of irreducible  complex
representations of $G$. Let $\nu:H\to\GL(\mfW)$ be a finite-dimensional
representation of $H$ and let $\Hom_H(\mfV_\sfa,\mfW)$ be the space of
all $H$-equivariant homomorphisms $L:\mfV_\sfa\to\mfW$. We embed
$\mfV_\sfa\otimes\Hom_H(\mfV_\sfa,\mfW)$ into the space
$\Sec(G\times_\nu\mfW)$ of smooth sections of the associated vector
bundle $G\times_\nu\mfW$ by assigning to
$\sfv\otimes L\in\mfV_\sfa\otimes\Hom_H(\mfV_\sfa,\mfW)$ the
$H$-equivariant map
\begin{displaymath}
  g\in G\mapsto L\left(\tau_\sfa\left(g^{-1}\right)\sfv\right)\in\mfW\;.
\end{displaymath}

For a proof of the theorem below, we refer to \cite{wallach}.

\begin{tm}[Frobenius reciprocity]
  \label{tm:frobenius}
  The space $\Sec(G\times_\nu\mfW)$ decomposes into the direct sum
  \begin{displaymath}
    \sum_{\sfa\in\sfA}\mfV_\sfa\otimes\Hom_H(\mfV_\sfa,\mfW)\;.
  \end{displaymath}
  \hspace*{\fill} $\square$
\end{tm}

Applying Theorem~\ref{tm:frobenius} to the
representation $\boldsymbol{u}_l$, we obtain a decomposition of
$\Sec(Q_l)$ into the sum
\begin{displaymath}
  \sum_{\sfa\in\sfA}\mfV_\sfa\otimes\Hom_H(\mfV_\sfa,\mfW_l)\;.
\end{displaymath}
Using this and Proposition~\ref{sz:splitting-Q}, we now want to find
invariant subspaces for the symplectic Dirac operators $D$ and $\tilde D$.
We start with expressing the symplectic Clifford multiplication
by means of the Hermite functions $\boldsymbol{h}_\alpha$.

\begin{lm}
  \label{lm:cliffm-ny-hermfunc}
  For any $\alpha\in\N_0^n$ and $j=1,\ldots,n$,
  \begin{displaymath}
    X_j\cdot\boldsymbol{h}_\alpha
    =-\ri\alpha_j\boldsymbol{h}_{\alpha-\langle j\rangle}
    -\frac{\ri}{2}\boldsymbol{h}_{\alpha+\langle j\rangle}
  \end{displaymath}
  and
  \begin{displaymath}
    X_{n+j}\cdot\boldsymbol{h}_\alpha
    =-\alpha_j\boldsymbol{h}_{\alpha-\langle j\rangle}
    +\frac{1}{2}\boldsymbol{h}_{\alpha+\langle j\rangle}\;,
  \end{displaymath}
  where $\langle j\rangle=(\delta_{1j},\ldots,\delta_{nj})$.
\end{lm}
\begin{beweis}
  This follows from Equations~(\ref{eq:def-cliff-mult}) and
  (\ref{eq:def-herm-func}) and the relations
  \begin{displaymath}
    \boldsymbol{h}_l'(t)-t\boldsymbol{h}_l(t)
    =\boldsymbol{h}_{l+1}(t)
    \quad\mbox{and}\quad
    \boldsymbol{h}_l'(t)+t\boldsymbol{h}_l(t)
    =-2l\boldsymbol{h}_{l-1}(t)
  \end{displaymath}
  for the classical Hermite functions $\boldsymbol{h}_l$.
\end{beweis}

\begin{lm}
  \label{lm:for-def-Da}
  Let $L\in\Hom_H(\mfV_\sfa,\mfW_l)$ and set
  \begin{equation}
    \label{eq:def-Da}
    \mcD_\sfa(L)=-\sum_{j=1}^n X_j\cdot L\circ(\tau_\sfa)_\ast(X_{n+j})
    +\sum_{j=1}^n X_{n+j}\cdot L\circ(\tau_\sfa)_\ast(X_j)\;.
  \end{equation}
  Then
  \begin{displaymath}
    \mcD_\sfa(L)\in\Hom_H(\mfV_\sfa,\mfW_{l-1})
    \oplus\Hom_H(\mfV_\sfa,\mfW_{l+1})
  \end{displaymath}
  with the convention that $\mfW_{-1}=0$.
\end{lm}
\begin{beweis}
  By Lemma~\ref{lm:cliffm-ny-hermfunc},
  \begin{displaymath}
    X\cdot\mfW_l\subset\mfW_{l-1}\oplus\mfW_{l+1}
  \end{displaymath}
  for $X\in\mfm$. Consequently, $\mcD_\sfa(L)$ is a homomorphism from
  $\mfV_\sfa$ to $\mfW_{l-1}\oplus\mfW_{l+1}$. It remains to show that
  $\mcD_\sfa(L)$ is $H$-equivariant. Let $h\in H$ and $X,Y\in\mfm$. Then
  \begin{align*}
    \tau_\sfa\left(h^{-1}\right)\circ(\tau_\sfa)_\ast(X)\circ\tau_\sfa(h)
    &=\left.\frac{\rd}{\rd t}\tau_\sfa\left(h^{-1}\right)
      \circ\tau_\sfa(\exp(tX))\circ\tau_\sfa(h)\right|_{t=0}\\
    &=\left.\tau_\sfa\left(\exp\left(t \Ad\left(h^{-1}\right)X\right)\right)
    \right|_{t=0}\\
    &=(\tau_\sfa)_\ast\left(\kappa\left(h^{-1}\right)X\right)\;.
  \end{align*}
  Hence
  \begin{equation}
    \label{eq:equiv-Dk-1}
    (\tau_\sfa)_\ast(X)\circ\tau_\sfa(h)
    =\tau_\sfa(h)\circ(\tau_\sfa)_\ast
    \left(\kappa\left(h^{-1}\right)X\right)\;.
  \end{equation}
  Since the symplectic Clifford multiplication is $\Mp(\mfm)$-equivariant,
  i.e.
  \begin{displaymath}
    \rho(a)X\cdot\boldsymbol{m}(a)f=\boldsymbol{m}(a)(X\cdot f)
  \end{displaymath}
  for $a\in\Mp(\mfm)$ and $f\in L^2(\R^n)$, we furthermore have
  \begin{equation}
    \label{eq:equiv-Dk-2}
    X\cdot\lambda(h)f
    =\lambda(h)\left(\kappa\left(h^{-1}\right)X\cdot f\right)\;.
  \end{equation}
  Equations~(\ref{eq:equiv-Dk-1}) and (\ref{eq:equiv-Dk-2}) imply
  \begin{align*}
    X\cdot L\circ(\tau_\sfa)_\ast(Y)\circ\tau_\sfa(h)
    &=X\cdot L\circ\tau_\sfa(h)
    \circ(\tau_\sfa)_\ast\left(\kappa\left(h^{-1}\right)Y\right)\\
    &=X\cdot\lambda(h)\circ L
    \circ(\tau_\sfa)_\ast\left(\kappa\left(h^{-1}\right)Y\right)\\
    &=\lambda(h)\circ\left(\kappa\left(h^{-1}\right)X\cdot L
      \circ(\tau_\sfa)_\ast\left(\kappa\left(h^{-1}\right)Y\right)\right)\;.
  \end{align*}
  Since the definition of $\mcD_\sfa(L)$ does not depend on the choice
  of the symplectic basis $(X_1,\ldots,X_{2n})$ and since
  $(\kappa(h)X_1,\ldots,\kappa(h)X_{2n})$ for any $h\in H$ is
  again a symplectic basis, this proves the lemma.
\end{beweis}

\medskip

For $\sfa\in\sfA$, set
\begin{displaymath}
  \bfU_\sfa=\sum_{l=0}^\infty\Hom_H(\mfV_\sfa,\mfW_l)\;.
\end{displaymath}
By Lemma~\ref{lm:for-def-Da}, Equation~(\ref{eq:def-Da}) defines a
homomorphism $\mcD_\sfa:\bfU_\sfa\to\bfU_\sfa$. The symplectic
Dirac operator $\mcD$ and the homomorphisms $\mcD_\sfa$ are related by

\begin{lm}
  \label{lm:rel-D-Da}
  If $\sfv\otimes L\in\mfV_\sfa\otimes\bfU_\sfa$, then
  \begin{displaymath}
    \mcD(\sfv\otimes L)=\sfv\otimes\mcD_\sfa(L)\;.
  \end{displaymath}
\end{lm}
\begin{beweis}
  We deduce that
  \begin{align*}
    X(\sfv\otimes L)(g)
    &=\left.\frac{\rd}{\rd t}(\sfv\otimes L)(g\exp(tX))\right|_{t=0}\\
    &=\left.\frac{\rd}{\rd t}
      L\left(\tau_\sfa\left(\exp(-tX)g^{-1}\right)\sfv\right)\right|_{t=0}\\
    &=-L\circ(\tau_\sfa)_\ast(X)
    \left(\tau_\sfa\left(g^{-1}\right)\sfv\right)
  \end{align*}
  for $X\in\mfm$ and $g\in G$. This together with
  Proposition~\ref{sz:D-on-ssymsp} gives the assertion.
\end{beweis}

\medskip

The decompositions obtained above and Lemma~\ref{lm:rel-D-Da} allow us
to consider the symplectic Dirac operator $\mcD$ as a map
\begin{displaymath}
  \mcD:\sum_{\sfa\in\sfA}\mfV_\sfa\otimes\bfU_\sfa
  \to\sum_{\sfa\in\sfA}\mfV_\sfa\otimes\bfU_\sfa\;.
\end{displaymath}
Moreover, we have

\begin{sz}
  \label{sz:inv-spaces-D}
  For any $\sfa\in\sfA$, the space $\mfV_\sfa\otimes\bfU_\sfa$ is
  invariant under $\mcD$.
\end{sz}
\begin{beweis}
  This is a consequence of Lemmas~\ref{lm:for-def-Da} and
  \ref{lm:rel-D-Da}.
\end{beweis}

\medskip

The same can be done for the symplectic Dirac operator $\tilde\mcD$ and the
second order operator $\mcP=\ri\left[\tilde\mcD,\mcD\right]$.

\begin{lm}
  \label{lm:for-def-tildeDa}
  Let $L\in\Hom_H(\mfV_\sfa,\mfW_l)$ and set
  \begin{equation}
    \label{eq:def-tildeDa}
    \tilde\mcD_\sfa(L)=-\sum_{j=1}^{2n} X_j
    \cdot L\circ(\tau_\sfa)_\ast(X_j)\;,
  \end{equation}
  where here the basis $(X_1,\ldots,X_{2n})$ of $\mfm$ is assumed to be
  unitary. Then
  \begin{displaymath}
    \tilde\mcD_\sfa(L)\in\Hom_H(\mfV_\sfa,\mfW_{l-1})
    \oplus\Hom_H(\mfV_\sfa,\mfW_{l+1})\;.
  \end{displaymath}
\end{lm}
\begin{beweis}
  Using that $\kappa(h)$ is unitary for all $h\in H$, the proof is the
  same as the proof of Lemma~\ref{lm:for-def-Da}.
\end{beweis}

\medskip

Consequently, Equation~(\ref{eq:def-tildeDa}) gives rise to a
homomorphism $\tilde\mcD_\sfa:\bfU_\sfa\to\bfU_\sfa$.

\begin{lm}
  \label{lm:rel-tildeD-tildeDa}
  Let $\sfv\otimes L\in\mfV_\sfa\otimes\bfU_\sfa$. Then
  \begin{displaymath}
    \tilde\mcD(\sfv\otimes L)=\sfv\otimes\tilde\mcD_\sfa(L)
    \quad\mbox{and}\quad
    \mcP(\sfv\otimes L)=\sfv\otimes\mcP_\sfa(L)\;,
  \end{displaymath}
  where $\mcP_\sfa=\ri\left[\tilde\mcD_\sfa,\mcD_\sfa\right]$.
\end{lm}
\begin{beweis}
  This follows from Proposition~\ref{sz:tildeD-on-ssymsp} and
  Lemma~\ref{lm:rel-D-Da}.
\end{beweis}

\begin{sz}
  The spaces $\mfV_\sfa\otimes\bfU_\sfa$ are invariant under
  $\tilde\mcD$ and $\mcP$.
\end{sz}
\begin{beweis}
  This is immediate from Lemmas~\ref{lm:for-def-tildeDa} and
  \ref{lm:rel-tildeD-tildeDa} and Proposition~\ref{sz:inv-spaces-D}.
\end{beweis}

\begin{bs}
  We consider the complex projective space $\C P^n=\SU(n+1)/\rU(n)$,
  where the group $\rU(n)$ is embedded into $\SU(n+1)$ by means of
  \begin{displaymath}
    B\in\rU(n)\mapsto
    \begin{pmatrix}
      \det(B)^{-1} & 0\\
      0 & B
    \end{pmatrix}
    \in\SU(n+1)\;.
  \end{displaymath}
  Let $\mfsu(n+1)$ and $\mfu(n)$ denote the Lie algebras of $\SU(n+1)$
  and $\rU(n)$, respectively, and let $\mfm$ be the image of the
  homomorphism
  \begin{displaymath}
    \Phi:z\in\C^n\mapsto
    \begin{pmatrix}
      0 & -\bar z^\rT\\
      z & 0
    \end{pmatrix}
    \in\mfsu(n+1)\;.
  \end{displaymath}
  Then $\mfm$ is an $\ad(\mfu(n))$-invariant complement of
  $\mfu(n)\subset\mfsu(n+1)$ satisfying $[\mfm,\mfm]\subset\mfu(n)$.
  We identify $\mfm$ with $\C^n$ via $\Phi$. With this, the
  isotropy representation $\kappa$ takes the form
  \begin{displaymath}
    \kappa(B)=\det(B)B
  \end{displaymath}
  for $B\in\rU(n)$. Moreover, $\omega_0$, $J_0$ and $\mathbf{g}_0$
  are the standard symplectic form, complex structure and inner product
  on $\C^n$, i.e.
  \begin{displaymath}
    \omega_0(z,w)=\operatorname{Im}\left(\bar z^\rT w\right)\;,\quad
    J_0 z=\ri z\quad\mbox{and}\quad
    \mathbf{g}_0(z,w)=\operatorname{Re}\left(\bar z^\rT w\right)
  \end{displaymath}
  for $z,w\in\C^n$. There exists a lift $\tilde\kappa$ of $\kappa$ iff
  $n$ is odd and in this case it is unique (for more details,
  see \cite{wyss}). Since $\kappa$ maps onto $\rU(n)$ and since
  the representations $\boldsymbol{u}_l$ are irreducible and pairwise
  inequivalent, by Schur's lemma, $\Hom_H(\mfV_\sfa,\mfW_l)$ is
  non-trivial for only finitely many $\sfa\in\sfA$. Therefore,
  the spaces $\mfV_\sfa\otimes\bfU_\sfa$ have finite dimension.
\end{bs}

In the next section, we examine the case of the complex projective line
$\C P^1=\SU(2)/\rU(1)$ in more detail.

\section{Symplectic Dirac operators on $\C P^1$
\label{sec:cp1}}

We set
\begin{displaymath}
  E_0=
  \begin{pmatrix}
    \ri & 0\\
    0 & \ri
  \end{pmatrix}
  \,,\quad
  E_1=
  \begin{pmatrix}
    0 & 1\\
    -1 & 0
  \end{pmatrix}
  \,,\quad
  E_2=
  \begin{pmatrix}
    0 & \ri\\
    \ri & 0
  \end{pmatrix}
  \,.
\end{displaymath}
The matrices $E_0,E_1,E_2$ form a basis of the Lie algebra $\mfsu(2)$
and satisfy
\begin{equation}
  \label{eq:kommu-su2}
  [E_0,E_1]=2 E_2\;,\quad
  [E_0,E_2]=-2 E_1\;,\quad
  [E_1,E_2]=2 E_0\;.
\end{equation}
Furthermore, the subalgebra $\mfu(1)\subset\mfsu(2)$ is spanned by $E_0$
and $\{E_1,E_2\}$ is a symplectic basis of the complement $\mfm$.
Since the complex structure $J_0$ on $\mfm$ is given by $J_0 E_1=E_2$,
this basis is also unitary.
We set $\jmath_0=\rho_\ast^{-1}(J_0)$ and define a homomorphism
$\tilde\kappa:\rU(1)\to\hat\rU(\mfm)$ by
\begin{displaymath}
  \tilde\kappa_\ast(E_0)=2\jmath_0\;.
\end{displaymath}
By Equation~(\ref{eq:kommu-su2}),
\begin{displaymath}
  \kappa_\ast(E_0)=
  \begin{pmatrix}
    0 & -2\\
    2 & 0
  \end{pmatrix}
  =2 J_0\;.
\end{displaymath}
Consequently, $\rho_\ast\circ\tilde\kappa_\ast=\kappa_\ast$, which implies
$\rho\circ\tilde\kappa=\kappa$. Hence $\tilde\kappa$ is a lift of $\kappa$.

Let $\mfV_k$ for $k\in\N_0$ be the vector space of homogeneous
polynomials of degree $k$ in two variables $z_1$ and $z_2$ with complex
coefficients. The polynomials $p_{k,0},\ldots,p_{k,k}$ defined by
\begin{displaymath}
  p_{k,j}(z_1,z_2)=z_1^{k-j}z_2^j
\end{displaymath}
form a basis of $\mfV_k$. Let $\tau_k:\SU(2)\to\GL(\mfV_k)$ be the
representation of $\SU(2)$ induced by the canonical action
of $\SU(2)$ on $\C^2$, i.e.
\begin{displaymath}
  (\tau_k(g)p)(z)=p\left(g^{-1}z\right)
\end{displaymath}
for $g\in\SU(2)$, $p\in\mfV_k$ and $z\in\C^2$.
As it is well-known (cf. e.g. \cite{zelobenko}), the representations
$\tau_k$ are irreducible and pairwise inequivalent
and these are all irreducible complex representations
of $\SU(2)$ up to isomorphism. Moreover, one has

\begin{lm}
  \label{lm:weights-tauk}
  For all $k\in\N_0$ and $j=0,\ldots,k$,
  \begin{align*}
    (\tau_k)_\ast(E_0)p_{k,j}&=\ri(2j-k)p_{k,j}\;,\\
    (\tau_k)_\ast(E_1)p_{k,j}&=(j-k)p_{k,j+1}+jp_{k,j-1}\;,\\
    (\tau_k)_\ast(E_2)p_{k,j}&=\ri(j-k)p_{k,j+1}-\ri jp_{k,j-1}\;.
  \end{align*}
  \hspace*{\fill} $\square$
\end{lm}

\begin{lm}
  \label{lm:calcul-invsp}
  Let $k,l\in\N_0$. If $(k+1)/2+l\in\{0,\ldots,k\}$, then the
  space $\Hom_{\rU(1)}(\mfV_k,\mfW_l)$ is one-dimensional and
  generated by the linear map $L_{k,l}:\mfV_k\to\mfW_l$ defined by
  \begin{displaymath}
    L_{k,l}(p_{k,j})=\left\{
      \begin{array}{cl}
        \boldsymbol{h}_l & \mbox{for}\quad j=(k+1)/2+l\\
        0 & \mbox{else}
      \end{array}
    \right..
  \end{displaymath}
  In all other cases, $\Hom_{\rU(1)}(\mfV_k,\mfW_l)$ is trivial.
\end{lm}
\begin{beweis}
  Let $\mathsf{H}_0:L^2(\R)\to L^2(\R)$ be the Hamilton operator of
  the one-dimensional harmonic oscillator, i.e.
  \begin{displaymath}
    (\mathsf{H}_0 f)(x)=\frac{1}{2}
    \left(\frac{\rd^2 f}{\rd x^2}(x)-x^2f(x)\right)\;.
  \end{displaymath}
  Then (cf. \cite{hab2-buch})
  \begin{displaymath}
    \boldsymbol{m}_\ast(\jmath_0)=-\ri\mathsf{H}_0
  \end{displaymath}
  and
  \begin{displaymath}
    \mathsf{H}_0\boldsymbol{h}_l=-\frac{2l+1}{2}\boldsymbol{h}_l\;.
  \end{displaymath}
  Hence
  \begin{displaymath}
    (\boldsymbol{u}_l\circ\tilde\kappa)_\ast(E_0)\boldsymbol{h}_l
    =2\boldsymbol{m}_\ast(\jmath_0)\boldsymbol{h}_l
    =-2\ri\mathsf{H}_0\boldsymbol{h}_l
    =\ri(2l+1)\boldsymbol{h}_l\;.
  \end{displaymath}
  Applying Schur's lemma, the assertion now follows from
  Lemma~\ref{lm:weights-tauk}.
\end{beweis}

\medskip

As shown in the previous section, the spaces $\mfV_k\otimes\bfU_k$ with
\begin{displaymath}
  \bfU_k=\sum_{l=0}^\infty\Hom_{\rU(1)}(\mfV_k,\mfW_l)
\end{displaymath}
are invariant under the symplectic Dirac operators $\mcD$ and
$\tilde\mcD$.

\begin{fo}
  If $k$ is odd, then
  \begin{displaymath}
    \bfU_k=\sum_{l=0}^{(k-1)/2}\Hom_{\rU(1)}(\mfV_k,\mfW_l)
  \end{displaymath}
  and $\dim(\mfV_k\otimes\bfU_k)=(k+1)^2/2$. If $k$ is even, then
  $\bfU_k$ is trivial.
\end{fo}
\begin{beweis}
  This is a consequence of Lemma~\ref{lm:calcul-invsp} and
  $\dim\mfV_k=k+1$.
\end{beweis}

\medskip

According to Lemmas~\ref{lm:rel-D-Da} and \ref{lm:rel-tildeD-tildeDa},
in order to calculate the symplectic Dirac operators $\mcD$ and
$\tilde\mcD$ on $\C P^1$, we have to compute the linear maps
$\mcD_k,\tilde\mcD_k:\bfU_k\to\bfU_k$ given by
\begin{displaymath}
  \mcD_k(L)=-E_1\cdot L\circ(\tau_k)_\ast(E_2)
  +E_2\cdot L\circ(\tau_k)_\ast(E_1)
\end{displaymath}
and
\begin{displaymath}
  \tilde\mcD_k(L)=-E_1\cdot L\circ(\tau_k)_\ast(E_1)
  -E_2\cdot L\circ(\tau_k)_\ast(E_2)\;.
\end{displaymath}
Thereby, we may assume that $k$ is odd. Then the maps
$L_{k,0},\ldots,L_{k,(k-1)/2}$ defined in Lemma~\ref{lm:calcul-invsp}
form a basis of $\bfU_k$.

\begin{sz}
  \label{sz:descr-Dk}
  For $l=0,\ldots,(k-1)/2$,
  \begin{displaymath}
    \mcD_k(L_{k,l})=l(k+1-2l)L_{k,l-1}
    +\left(\frac{k+1}{2}+l+1\right)L_{k,l+1}
  \end{displaymath}
  and
  \begin{displaymath}
    \tilde\mcD_k(L_{k,l})=-\ri l(k+1-2l)L_{k,l-1}
    +\ri\left(\frac{k+1}{2}+l+1\right)L_{k,l+1}\;.
  \end{displaymath}
\end{sz}
\begin{beweis}
  Assume that $j=(k+1)/2+l$. Then, by Lemmas~\ref{lm:cliffm-ny-hermfunc}
  and \ref{lm:weights-tauk}, we get
  \begin{align*}
    \mcD_k(L_{k,l})(p_{k,j-1})
    &=-E_1\cdot L_{k,l}((\tau_k)_\ast(E_2)p_{k,j-1})
    +E_2\cdot L_{k,l}((\tau_k)_\ast(E_1)p_{k,j-1})\\
    &=-E_1\cdot L_{k,l}(\ri(j-k-1)p_{k,j}-\ri(j-1)p_{k,j-2})\\
    &\hspace{5mm}
    +E_2\cdot L_{k,l}((j-k-1)p_{k,j}+(j-1)p_{k,j-2})\\
    &=(j-k-1)(-\ri E_1\cdot\boldsymbol{h}_l
    +E_2\cdot\boldsymbol{h}_l)\\
    &=-2l(j-k-1)\boldsymbol{h}_{l-1}\\
    &=l(k+1-2l)\boldsymbol{h}_{l-1}
  \end{align*}
  and
  \begin{align*}
    \mcD_k(L_{k,l})(p_{k,j+1})
    &=-E_1\cdot L_{k,l}((\tau_k)_\ast(E_2)p_{k,j+1})
    +E_2\cdot L_{k,l}((\tau_k)_\ast(E_1)p_{k,j+1})\\
    &=-E_1\cdot L_{k,l}(\ri(j-k+1)p_{k,j+2}-\ri(j+1)p_{k,j})\\
    &\hspace{5mm}
    +E_2\cdot L_{k,l}((j-k+1)p_{k,j+2}+(j+1)p_{k,j})\\
    &=(j+1)(\ri E_1\cdot\boldsymbol{h}_l
    +E_2\cdot\boldsymbol{h}_l)\\
    &=(j+1)\boldsymbol{h}_{l+1}\\
    &=\left(\frac{k+1}{2}+l+1\right)\boldsymbol{h}_{l+1}\;.
  \end{align*}
  Analogously, we see
  \begin{displaymath}
    \tilde\mcD_k(L_{k,l})(p_{k,j-1})
    =-\ri l(k+1-2l)\boldsymbol{h}_{l-1}
  \end{displaymath}
  and
  \begin{displaymath}
    \tilde\mcD_k(L_{k,l})(p_{k,j+1})
    =\ri\left(\frac{k+1}{2}+l+1\right)\boldsymbol{h}_{l+1}\;.
  \end{displaymath}
  This together with Lemmas~\ref{lm:for-def-Da} and
  \ref{lm:for-def-tildeDa} yields the claim.
\end{beweis}

\medskip

We now normalize the homomorphisms $L_{k,0},\ldots,L_{k,(k-1)/2}$ by
\begin{displaymath}
  L_{k,l}^\circ:=\sqrt{\frac{((k+1)/2+l)!\,((k-1)/2-l)!}{2^l l!}}L_{k,l}
\end{displaymath}
and set
\begin{displaymath}
  a_{k,l}:=\sqrt{2l\left(\frac{(k+1)^2}{4}-l^2\right)}\;.
\end{displaymath}
Then a straightforward calculation gives

\begin{fo}
  \label{fo:descr-Dk-herm}
  For $l=0,\ldots,(k-1)/2$,
  \begin{displaymath}
    \mcD_k(L_{k,l}^\circ)
    =a_{k,l}L_{k,l-1}^\circ+a_{k,l+1}L_{k,l+1}^\circ
  \end{displaymath}
  and
  \begin{displaymath}
    \tilde\mcD_k(L_{k,l}^\circ)
    =-\ri a_{k,l}L_{k,l-1}^\circ+\ri a_{k,l+1}L_{k,l+1}^\circ\;.
  \end{displaymath}
  In particular, the matrix representations of $\mcD_k$ and $\tilde\mcD_k$
  with respect to the basis
  $\left\{L_{k,0}^\circ,\ldots,L_{k,(k-1)/2}^\circ\right\}$
  are Hermitian.
  \hspace*{\fill} $\square$
\end{fo}

We draw the following conclusions from Corollary~\ref{fo:descr-Dk-herm}.
First, for the operators $\mcP_k=\ri\left[\tilde\mcD_k,\mcD_k\right]$,
we obtain

\begin{fo}
  For $l=0,\ldots,(k-1)/2$,
  \begin{align*}
    \mcP_k(L_{k,l}^\circ)
    &=2\left(a_{k,l+1}^2-a_{k,l}^2\right)L_{k,l}^\circ\\
    &=\left((k+1)^2-3(2l+1)^2-1\right)L_{k,l}^\circ\;.
  \end{align*}
  \hspace*{\fill} $\square$
\end{fo}

\begin{fo}
  \begin{enumerate}
  \item[{\rm (i)}] The operators $\mcD_k$ and $\tilde\mcD_k$ have the
    same eigenvalues.
  \item[{\rm (ii)}] The spectrum of $\mcD_k$ and $\tilde\mcD_k$ is
    symmetric, i.e. if $\lambda$ is an eigenvalue, so is $-\lambda$.
  \item[{\rm (iii)}] The kernels of $\mcD_k$ and $\tilde\mcD_k$ are
    one-dimensional, if $(k+1)/2$ is odd, and trivial otherwise.
  \end{enumerate}
\end{fo}
\begin{beweis}
  One easily verifies that the characteristic polynomials of $\mcD_k$
  and $\tilde\mcD_k$ are the same and of the form
  $\lambda p_k\left(\lambda^2\right)$, if $(k+1)/2$ is odd, and
  $p_k\left(\lambda^2\right)$ otherwise, for some polynomial $p_k$.
  Since $a_{k,1},\ldots,a_{k,(k-1)/2}$ are non-zero, the dimension of
  the kernel
  of $\mcD_k$ is at most $1$. Moreover, if $(k+1)/2$ is even, then
  \begin{displaymath}
    \det(\mcD_k)=\prod_{r=1}^{(k+1)/4}a_{k,2r-1}^2\;.
  \end{displaymath}
\end{beweis}

\begin{fo}
  The spectrum of the symplectic Dirac operators $\mcD$ and $\tilde\mcD$
  on $\C P^1$ is unbounded above and below.
\end{fo}
\begin{beweis}
  Let $\bfU_k$ be endowed with the norm defined by
  \begin{displaymath}
    \left\|\sum_{l=0}^{(k-1)/2}z_l L_{k,l}^\circ\right\|^2
    =\sum_{l=0}^{(k-1)/2}|z_l|^2
  \end{displaymath}
  and let $\|\mcD_k\|$ be the corresponding operator norm of $\mcD_k$.
  Since the matrix representation of $\mcD_k$ with respect to
  $\left\{L_{k,0}^\circ,\ldots,L_{k,(k-1)/2}^\circ\right\}$ is Hermitian
  and the spectrum of $\mcD_k$ is symmetric, $\|\mcD_k\|$ is the
  absolut value of the largest and smallest eigenvalue of $\mcD_k$.
  Now the assertion follows from
  \begin{displaymath}
    \|\mcD_k\|\ge\left\|\mcD_k L_{k,0}^\circ\right\|
    =a_{k,1}\ge\frac{k-1}{2}\;.
  \end{displaymath}
  
\end{beweis}

\begin{bm0}
  The spectrum of the second order operator $\mcP$ on $\C P^1$ was
  already computed in \cite{khab}. For the general case, the coincidence
  of the spectra of $\mcD$ and $\tilde\mcD$ and their symmetry was proven
  in \cite{hab-klein} (cf. also \cite{hab2-buch}) by means of a Fourier
  transform for symplectic spinor fields. For an explicit calculation of
  the eigenvalues and eigenvectors of $\mcD_k$ and $\tilde\mcD_k$
  for small $k$, see \cite{steffen}.
\end{bm0}

\end{document}